\PassOptionsToPackage{unicode}{hyperref}
\PassOptionsToPackage{hyphens}{url}
\documentclass[
  letterpaper,
  10pt,
  conference]{ieeeconf}
\usepackage[]{times}
\usepackage{amssymb,amsmath}
\usepackage{ifxetex,ifluatex}
\ifnum 0\ifxetex 1\fi\ifluatex 1\fi=0 
  \usepackage[T1]{fontenc}
  \usepackage[utf8]{inputenc}
  \usepackage{textcomp} 
\else 
  \usepackage{unicode-math}
  \defaultfontfeatures{Scale=MatchLowercase}
  \defaultfontfeatures[\rmfamily]{Ligatures=TeX,Scale=1}
\fi
\IfFileExists{upquote.sty}{\usepackage{upquote}}{}
\IfFileExists{microtype.sty}{
  \usepackage[]{microtype}
  \UseMicrotypeSet[protrusion]{basicmath} 
}{}
\makeatletter
\@ifundefined{KOMAClassName}{
  \IfFileExists{parskip.sty}{%
    \usepackage{parskip}
  }{
    \setlength{\parindent}{0pt}
    \setlength{\parskip}{6pt plus 2pt minus 1pt}}
}{
  \KOMAoptions{parskip=half}}
\makeatother
\usepackage{xcolor}
\IfFileExists{xurl.sty}{\usepackage{xurl}}{} 
\IfFileExists{bookmark.sty}{\usepackage{bookmark}}{\usepackage{hyperref}}
\hypersetup{
  hidelinks,
  pdfcreator={LaTeX via pandoc}}
\usepackage{xparse}
\usepackage{amsmath}
\usepackage{bm}
\usepackage{mathtools}

\usepackage[sort,noadjust]{cite}
\usepackage{color,soul}
\sethlcolor{cyan}
\usepackage{tabularx}
\IEEEoverridecommandlockouts
\usepackage{cases}
\usepackage{diagbox}
\usepackage{algorithm}
\usepackage{algcompatible}

\usepackage{url}
\usepackage{booktabs}
\usepackage{multirow}
\makeatletter
\@ifpackageloaded{subfig}{}{\usepackage{subfig}}
\@ifpackageloaded{caption}{}{\usepackage{caption}}
\AtBeginDocument{%

}
\AtBeginDocument{%

}
\@ifpackageloaded{float}{}{\usepackage{float}}
\floatstyle{ruled}
\@ifundefined{c@chapter}{\newfloat{codelisting}{h}{lop}}{\newfloat{codelisting}{h}{lop}[chapter]}
\floatname{codelisting}{Listing}

\makeatother

\title{\textbf{Optimistic planning for the near-optimal control of nonlinear switched discrete-time systems with stability guarantees}\\}
\author{Mathieu
Granzotto\thanks{Mathieu Granzotto, Romain Postoyan and Jamal Daafouz are with the Universit\'e de Lorraine, CNRS, CRAN, F-54000 Nancy, France (e-mails: \{name.surname\}@univ-lorraine.fr).},
Romain Postoyan, Lucian
Buşoniu\thanks{Lucian Bu\c{s}oniu is with the Department of Automation, Technical University of Cluj-Napoca, Memorandumului 28, 400114 Cluj-Napoca, Romania (e-mail: lucian.busoniu@aut.utcluj.ro). This work was supported by a grant of Romanian Ministry of Research and Innovation, CNCS - UEFISCDI, project number PN-III-P1-1.1-TE-2016-0670, within PNCDI III.},
Dragan
Nešić\thanks{Dragan Nešić is with the Department of Electrical and Electronic Engineering, University of Melbourne, Parkville, VIC 3010, Australia (e-mail: dnesic@unimelb.edu.au). His work was supported by the Australian Research Council under the Discovery Project DP170104099.},
and Jamal Daafouz}
\date{}

\begin{document}
\maketitle

\newcommand\numberthis{\addtocounter{equation}{1}\tag{\theequation}}

\def \proof{\begin{IEEEproof}}
\def \endproof{\end{IEEEproof}}

\def \wDelta{\widetilde\Delta}
\def \wdelta{\widetilde\delta}

\newcommand{\awb}{\overline{\alpha}_W}
\newcommand{\avb}{\overline{\alpha}_V}
\newcommand{\ayb}{\overline{\alpha}_Y}
\newcommand{\aybi}{\ayb^{-1}}
\newcommand{\bya}{\underline{\alpha}_Y}
\newcommand{\byai}{\bya^{-1}}
\newcommand{\aw }{\alpha_W}
\newcommand{\awi}{\aw^{-1}}
\newcommand{\ay }{{\alpha}_Y}
\newcommand{\ayi}{\ay^{-1}}
\newcommand{\ayw}{\widetilde{\alpha}_Y}
\newcommand{\ayh}{\widehat{\alpha}_Y}

\newcommand{\Reals}{\mathbb{R}}
\newcommand{\PReals}{\mathbb{R}_{\geq 0}}
\newcommand{\Uset}{\mathcal{U}}

\newcommand{\Positives}{\mathbb{Z}_{>0}}
\newcommand{\ZeroPositiv}{\mathbb{Z}_{\geq 0}}
\newcommand{\Wset}{\mathcal{W}}

\newcommand{\useq}{\bm{u}}
\NewDocumentCommand{\ustar}{o}{%
  \IfNoValueTF{#1}
    {\bm{u}_{\gamma,d}^{\pmb{*}}}
    {\bm{u}_{\gamma,#1}^{\pmb{*}}}
}

\NewDocumentCommand{\us}{o}{%
  \IfNoValueTF{#1}
    {u^*}
    {u^* _ {#1}}
}
\newcommand{\Usetstar}{\Uset^*_{\gamma,d}}

\newcommand{\uhat}{\bm{\hat{u}}_{\gamma,d}}
\newcommand{\Uhatset}{\widehat{\mathcal{U}}_{\gamma,d}}

\NewDocumentCommand{\phis}{o}{%
  \IfNoValueTF{#1}
    {\phi}
    {\phi^* _ {#1}}
}
\NewDocumentCommand{\ells}{o}{%
  \IfNoValueTF{#1}
    {\ell}
    {\ell^*_{#1}}
}

\newcommand{\K}{\mathcal{K}}
\newcommand{\KK}{\mathcal{KK}}
\newcommand{\KL}{\mathcal{KL}}

\newcommand{\I}{\mathbb{I}}

\NewDocumentCommand{\V}{o}{%
  \IfNoValueTF{#1}
    {V_{\gamma,d}}
    {V_{\gamma,#1}}%
}
\NewDocumentCommand{\J}{o}{%
  \IfNoValueTF{#1}
    {J_{\gamma,d}}
    {J_{\gamma,#1}}%
}
\NewDocumentCommand{\Y}{o}{%
  \IfNoValueTF{#1}
    {Y_{\gamma,d}}
    {Y_{\gamma,#1}}%
}

\NewDocumentCommand{\Vhat}{o}{%
  \IfNoValueTF{#1}
    {{\widehat V}_{\gamma,d}}
    {{\widehat V}_{\gamma,#1}}%
}

\NewDocumentCommand{\Yhat}{o}{%
  \IfNoValueTF{#1}
    {{\widehat Y}_{\gamma,d}}
    {{\widehat Y}_{\gamma,#1}}%
}

\newcommand{\gfactor}{\frac{1-\gamma}{\gamma}}
\newcommand{\tgfactor}{\tfrac{1-\gamma}{\gamma}}

\newcommand{\gdfactor}{\gfactor\frac{1}{1-\gamma^{d}}}
\newcommand{\tgdfactor}{\tgfactor\tfrac{1}{1-\gamma^{d}}}

\newcommand{\gstarfactor}{\frac{1-\gamma^*}{\gamma^*}}
\newcommand{\tgstarfactor}{\tfrac{1-\gamma^*}{\gamma^*}}

\newcommand{\gdstarfactor}{\gstarfactor\frac{1}{1-(\gamma^*)^{d^*}}}
\newcommand{\tgdstarfactor}{\tgstarfactor\tfrac{1}{1-(\gamma^*)^{d^*}}}

\DeclarePairedDelimiter\floor{\lfloor}{\rfloor}

\NewDocumentCommand{\uref}{m}{%
    {\textsubscript{(\ref{#1})}}
}

\newtheorem{SA}{Standing Assumption (SA)}
\renewcommand\theSA{\unskip}
\newtheorem{corollary}{Corollary}
\newtheorem{assumption}{Assumption}
\newtheorem{remark}{Remark}
\newtheorem{lemma}{Lemma}
\newtheorem{theorem}{Theorem}
\newtheorem{proposition}{Proposition}

\newtheorem{example}{Example}

\newtheorem{dremark}{\hl{Draft Remark}}

\begin{abstract}
Originating in the artificial intelligence literature, optimistic planning (OP) is an algorithm that generates near-optimal control inputs for generic nonlinear discrete-time systems whose input set is finite. This technique is therefore relevant for the near-optimal control of nonlinear switched systems, for which the switching signal is the control. However, OP exhibits several limitations, which prevent its application in a standard control context. First, it requires the stage cost to take values in [0,1], an unnatural prerequisite as it excludes,  for instance, quadratic stage costs.   Second, it requires the cost function to be discounted. Third, it applies for reward maximization, and not cost minimization. In this paper, we modify OP to overcome these limitations, and we call the new algorithm OPmin. We then make stabilizability and detectability assumptions, under which we derive near-optimality guarantees for OPmin and we show that the obtained bound has major advantages compared to the bound originally given by OP. In addition, we  prove that a system whose inputs are generated by OPmin in a receding-horizon fashion exhibits stability properties. As a result, OPmin provides a new tool for the near-optimal, stable control of nonlinear switched discrete-time systems for generic cost functions.

\end{abstract}

\hypertarget{introduction}{%
\section{Introduction}\label{introduction}}

Optimistic planning (OP) is an algorithm that computes near-optimal
control inputs for generic nonlinear discrete-time systems and
infinite-horizon discounted costs, provided the set of inputs is finite,
see \cite{hren2008,munos2014bandits}. Given the current state, OP
intelligently develops the tree of possible future states which are
enumerable, as the input set is finite. By prioritizing branches with
smaller costs, which are optimistic candidates to the infinite-horizon
cost, OP efficiently exploits the available computational power. It then
returns an optimal sequence of inputs for a \emph{finite}-horizon
discounted cost, where the horizon depends on the given computational
budget and on the state. Guarantees on the mismatch between the obtained
cost and the original infinite-horizon cost are provided in
\cite{hren2008} and are of the form \(\frac{\gamma^{d(x)}}{1-\gamma}\),
where \(\gamma\in(0,1)\) is the discount factor and \(d(x)\) is the
state-dependent horizon, which is related to the computation budget
\(B\).

OP is a priori well-suited for nonlinear switched discrete-time systems
for which the control input corresponds to the switching signal
\cite{BUSONIU2017297}. While the (near)-optimal control of switched
linear discrete-time systems is addressed in, e.g.,
\cite{abateDSLQR,5288566,6656817,DEAECTO20181,7588142}, the case of
nonlinear switched systems is still unraveling and concentrates on
continuous-time systems, see e.g.~\cite{Antsaklis2003,Zhu2015}. Even so,
algorithms are often presented for a particular class of systems,
consider finite-horizon optimality and ignore stability. There is
therefore a need for tools for the (near-)optimal control of nonlinear
switched systems. We propose a solution based on OP in this paper.

It appears that we cannot apply OP ``off-the-shelf'' adequately for
optimal control problems. Indeed, OP requires that: (i) the stage cost
takes value in \([0,1]\), which is not natural in control, as this
excludes quadratic stage costs, for instance; (ii) the cost is
discounted; (iii) the goal is to maximize the value function, and
adapting OP to minimization is not straightforward. We therefore modify
OP to overcome these limitations. We call this new algorithm OPmin.
Similar to OP, OPmin returns a sequence of inputs, which minimizes a
finite-horizon cost more efficiently than a brute-force approach (in
general).

We make stabilizability and detectability assumptions, based on which we
analyze the near-optimality guarantees of OPmin, that is, how the
computed finite-horizon cost function compares to the infinite-horizon
cost. The obtained bound on the mismatch between the two costs have the
next desirable features: (i) it does not explode for \(\gamma=1\),
contrary to the bound in \cite{hren2008}; (ii) it decreases as the state
is close to a given attractor, while the bound
\(\frac{\gamma^{d(x)}}{1-\gamma}\) in \cite{hren2008} is a constant for
a constant horizon. In addition, inspired by our recent work
\cite{granzotto2018,granzotto2019}, we address the question of
stability, which is ignored in \cite{BUSONIU2017297,hren2008}. For this
purpose, we rely on the same stabilizability and detectability
assumptions as for the near-optimality analysis. We prove that a system,
for which the inputs are generated by OPmin in a receding-horizon
fashion, satisfies a semiglobal practical stability property, where the
adjustable parameters are the computational budget of OPmin and the
possible discount factor. We use a generic measuring function to define
stability as in \cite{grimm2005,romain2016,granzotto2019}, thus covering
point and set stability in a unified way. By strengthening the
assumptions, we also derive a global exponential stability property.
These stability results differ from our recent works in
\cite{granzotto2018,granzotto2019} as the horizon here is
state-dependent, and not fixed, like in
\cite{granzotto2018,granzotto2019}.

Finally, we investigate the relationship between the original
infinite-horizon optimal value function, and the actual cost function
obtained by applying OPmin in a receding-horizon fashion to the system,
also known as running cost \cite{gruneperformance}. Assuming that the
closed-loop system satisfies a global exponential stability property
(for which we provide sufficient conditions, as mentioned above), we
show that, indeed, increasing the horizon of optimization, i.e.~the
computational budget, implies that the running cost approaches the
infinite-horizon optimal value function. Moreover, the mismatch between
the running cost and the infinite-horizon cost decreases exponentially
under mild conditions. We also provide the relative performance of the
running cost to the infinite-horizon value function, similarly explored
in \cite{gruneperformance}. Contrary to \cite{gruneperformance}, we do
not rely on relaxed dynamic programming assumptions for this purpose but
on the aforementioned general stabilizability and detectability
conditions. An example is provided to illustrate the theoretical
results.

We think that this paper conveys an important message. It illustrates
how an optimal algorithm from a different research field, namely
artificial intelligence, can be adapted and tailored to solve an
important control problem, here the near-optimal control of nonlinear
switched discrete-time systems. It also demonstrates how control
requirements, like stabilizability, detectability and stability, can be
exploited to improve the original near-optimality guarantees of the
algorithm.

It must be noted that tree-based algorithms have been considered in the
literature for switched systems, albeit with different purposes. In
\cite{dellarossatree}, the stability of linear switched systems under
arbitrarily switching is investigated for instance. The work in
\cite{1039806} considers a branch-and-bound approach for the
discrete-time optimal control of switched linear systems and quadratic
costs. On the other hand, (relaxed) dynamic approaches were considered
in \cite{1687497,4982639}. In particular, \cite{1687497} approximates
the infinite-horizon optimal control problem for linear switched
systems, and \cite{4982639} develops a value iteration approach
exploiting homogeneity of the system and stage costs. The main
difference between our present paper and these references is that we
address nonlinear switched systems and generic (discounted) costs.

The rest of the paper is organized as follows. Section
\ref{problem-statement} formally states the problem. OPmin is presented
in Section \ref{opmin}, and its near-optimality and stability properties
are analyzed in Section \ref{main-results}. Section \ref{example}
provides an example.

\(\\\) \noindent\textbf{Notation.} Let
\(\mathbb{R}:= (-\infty,\infty)\), \(\mathbb{R}_{\geq 0}:= [0,\infty)\),
\(\mathbb{Z}_{\geq 0}:= \{0,1,2,\ldots\}\) and
\(\mathbb{Z}_{>0}:= \{1,2,\ldots\}\). We use \((x,y)\) to denote
\([x^T,y^T]^T\), where \((x,y) \in \mathbb{R}^n\times\mathbb{R}^m\) and
\(n,m\in\mathbb{Z}_{>0}\). A function
\(\chi : \mathbb{R}_{\geq 0}\to \mathbb{R}_{\geq 0}\) is of class
\(\mathcal{K}\) if it is continuous, zero at zero and strictly
increasing, and it is of class \(\mathcal{K}_\infty\) if it is of class
\(\mathcal{K}\) and unbounded. A continuous function
\(\beta: \mathbb{R}_{\geq 0}\times\mathbb{R}_{\geq 0}\to \mathbb{R}_{\geq 0}\)
is of class \(\mathcal{KL}\) when \(\beta(\cdot,t)\) is of class
\(\mathcal{K}\) for any \(t\geq0\) and \(\beta(s,\cdot)\) is decreasing
to 0 for any \(s\geq0\). The notation \(\mathbb{I}\) stands for the
identity map from \(\mathbb{R}_{\geq 0}\) to \(\mathbb{R}_{\geq 0}\).
For any sequence \(\bm{u}=[u_0,u_1,\dots]\) of length
\(d\in\mathbb{Z}_{\geq 0}\cup\{\infty\}\) where
\(u_i \in \mathbb{R}^m\), \(i \in \{0,\ldots,d\}\), and any
\(k\in\{0,\ldots,d\}\), we use \(\bm{u}|_k\) to denote the first \(k\)
elements of \(\bm{u}\), i.e.~\(\bm{u}|_k = [u_0,\dots,u_{k-1}]\) and
\(\bm{u}|_0=\varnothing\) by convention. Let \(f\  :\
\mathbb{R}\to\mathbb{R}\), we use \(f^{(k)}\) for the composition of
function \(f\) to itself \(k\) times, where \(k\in\mathbb{Z}_{\geq 0}\),
and \(f^{(0)}=\mathbb{I}\). The Euclidean norm of a vector
\(x\in\mathbb{R}^n\) is denoted by \(|x|\). The distance of a vector
\(x\in\mathbb{R}^n\) to set \(\cal A\) is defined as
\(|x|_{\cal A} =\inf\{|z-x| : z\in\cal A\}\).

\hypertarget{problem-statement}{%
\section{Problem Statement}\label{problem-statement}}

Consider the system
\begin{equation}x_{k+1} = f_{u_{k}}(x_{k}),\label{eq:sys}\end{equation}
with state \(x \in \mathbb{R}^n\), input \(u \in \mathcal{U}\), where
\(\mathcal{U}:=\{1,\ldots,M\}\) is a \emph{finite} set of admissible
inputs with \(M\geq 2\), and \(f_u: \mathbb{R}^n \to \mathbb{R}^n\) for
every input \(u \in \mathcal{U}\). We use \(\phi(k,x,\bm{u}|_ k)\) to
denote the solution to system (\ref{eq:sys}) at time
\(k \in \mathbb{Z}_{\geq 0}\) with initial condition \(x\) and inputs
\(\bm{u}|_k=[u_{0},u_{1},\ldots,u_{k-1}]\), with the convention
\(\phi(0,x,\cdot)=\phi(0,x,\varnothing)=x\).

Our objective is to minimize the infinite-horizon cost \begin{equation}
\J[\infty](x,\bm{u}):=\sum\limits_{k=0}^{\infty}\gamma^{k} \ell_{u_k}(\phi(k,x,\bm{u}|_k)),
\label{eq:Jinfty}
\end{equation} where \(x\in\mathbb{R}^n\) is the initial state,
\(\bm{u}\) is an infinite sequence of admissible inputs,
\(\ell_u : \mathbb{R}^n \to \mathbb{R}_{\geq 0}\) is the stage cost
related to input \(u\in\mathcal{U}\), and \(\gamma \in (0,1]\) is the
discount factor, which may be equal to \(1\). Finding an infinite
sequence of inputs which minimizes (\ref{eq:Jinfty}) is very difficult
in general, as the case of linear switched systems with quadratic stage
cost already shows \cite{abateDSLQR}. We therefore aim at generating
sequences of inputs that approximately minimize (\ref{eq:Jinfty})
instead, in a sense made precise in the following. For this purpose, we
adapt optimistic planning (OP) as originally developed in
\cite{hren2008} to be applicable for: (i) stage costs which are not
constrained to take values in \([0,1]\), to cope with quadratic stage
costs for instance; (ii) the undiscounted case, i.e.~when \(\gamma=1\)
in (\ref{eq:Jinfty}); (iii) the minimization of (2), as opposed to
maximization. We call this new algorithm OPmin. Furthermore, we also aim
at ensuring stability properties for the induced closed-loop system. The
algorithm is presented in the next section.

\hypertarget{opmin}{%
\section{OPmin}\label{opmin}}

\hypertarget{main-idea}{%
\subsection{Main idea}\label{main-idea}}

The algorithm we are going to present minimizes exactly, as we will
prove, the following \emph{finite}-horizon cost \begin{equation}
  \J[d(x)](x,\bm{u}):=\sum\limits_{k=0}^{d(x)}\gamma^{k} \ell_{u_k}(\phi(k,x,\bm{u}|_k)),
\label{eq:J}\end{equation} where \(x\in\mathbb{R}^n\) is a given state,
\(d(x)\in\mathbb{Z}_{>0}\) is the horizon, which depends on \(x\), and
\(\bm{u}=[u_0,u_1,...,u_{d(x)}]\) is a \(d(x)+1\) sequence of admissible
inputs. The associated optimal value function is \begin{equation}
        \V[d(x)](x) := \min_{\bm{u}} \J[d(x)](x,\bm{u}),
\label{eq:V}\end{equation} and we define \(\bm{u}_{\gamma,d(x)}^*\) an
associated \emph{optimal input sequence},
i.e.~\(\V[d(x)](x)=\J[d(x)](x,\bm{u}_{\gamma,d(x)}^*(x))\).

Compared to minimizing (\ref{eq:Jinfty}), problem (\ref{eq:V}) with
finite \(d(x)\) is solvable, as the input set \(\mathcal{U}\) is finite.
A brute-force approach can do it by developing all possible sequences.
However, this is computationally intensive, in particular when \(d(x)\)
is large, as the computational cost grows exponentially with the
horizon. OPmin, on the other hand, can intelligently explore the
possible sequences to solve (\ref{eq:V}) with potentially larger
\(d(x)\) with the same computation, compared to a brute-force approach
\cite{hren2008}. As we will show next, longer horizons imply smaller
near-optimality bounds, and are therefore desirable.

\hypertarget{algorithm}{%
\subsection{Algorithm}\label{algorithm}}

The objective of the algorithm is to find an input sequence such that
(\ref{eq:J}) is minimized and \(\V[d(x)](x)\) `approximates well'
\(\V[\infty](x)\), as formalized later in Section \ref{main-results}. It
does so by exploring the possible choices of inputs
\emph{optimistically} until the exhaustion of given computational
resources. The computational resources available are denoted as a budget
\(B\), which corresponds to \(B+1\) `leaf expansions' (the root is
expanded even at \(B=0\)). We denote by \(\mathcal{T}\) the exploration
tree from initial state \(x\in\mathbb{R}^n\), constructed from
admissible input sequences and their respective cost. A leaf is a node
of \(\mathcal{T}\) with no children, and the the set of all leaves of
\(\mathcal{T}\) is denoted \(\mathcal{L}(\mathcal{T})\). At iteration
\(i\in\mathbb{Z}_{\geq 0}\), a leaf \(L_i\in\mathcal{L}(\mathcal{T})\)
is fully expanded. That is, for every \(u\in\mathcal{U}\), we add a
child to \(L_i\) labeled by the resulting state \(f_u(L_i)\). We denote
with a slight abuse of notation \(\bm{u}(L_i)\) the input sequence from
the root to the state of leaf \(L_i\). We denote by
\(J(L_i):=\J[d(i)](x,\bm{u}(L_i))\) cost (\ref{eq:J}) of the sequence
that takes \(x\) to the state of leaf \(L_i\), with
\(d(i)=\text{depth}(L_i)-1\), where depth is the number of edges (or
inputs) from the root to \(L_i\). The optimistic choice of leaf \(L_i\)
to expand is the leaf with minimal associated cost \(J\) of all
non-expanded leafs of \(\mathcal{T}\). The algorithm is formalized next.

\begin{algorithm}[H]
 \caption{\label{algo}Algorithm for OPmin}
 \begin{algorithmic}[1]
 \renewcommand{\algorithmicrequire}{\textbf{Input:}}
 \renewcommand{\algorithmicensure}{\textbf{Output:}}
 \REQUIRE budget $B$
 \ENSURE  depth explored $d(x)$, sequence $\bm{u}^*_{\gamma,d(x)}(x)$,  cost $\V[d(x)](x)$
 \STATEx \textit{Initialisation} :
 \STATE $d\leftarrow -1$
 \STATE tree  $\mathcal{T}\leftarrow \{[\,],0\}$ \hfill \COMMENT{the empty sequence and cost 0}
  \textit{Optimistic exploration}
  \FOR {$i = 0$ to  $B$ }
  \STATE find optimistic leaf $L_i \in \mathop\mathrm{arg\,min}\limits_{L\in\mathcal{L}(\mathcal{T})}J(L)$
 \STATEx \textit{add to $\mathcal{T}$ the children of $L_i$:}
  \STATE for each child  $c$ of $L_i$, $\mathcal{T}\leftarrow\mathcal{T}\cup\{\bm{u}(c),J(c)\}$
  \IF {$d<\text{depth}(L_i)-1$}
 \STATEx \textit{Leaf selection}
  \STATE  $S\leftarrow L_i$
  \STATE  $d\leftarrow \text{depth}(L_i)-1$
  \ENDIF
  \ENDFOR
 \STATE \textbf{return} $d(x)\leftarrow d$ and $\{\bm{u}_{\gamma,d(x)}^*(x),\V[d(x)](x)\}\leftarrow S$
 \end{algorithmic}
\end{algorithm}

The most notable steps of Algorithm \ref{algo} are lines 4-5, where the
optimistic exploration is realized. This optimistic choice guarantees
that any sequence from descendants of a node \(N\) will have costs \(J\)
greater than \(N\), as \(\ell\geq0\). This implies that the first leaf
to be expanded at a depth \(d'+1\) will be a suitable candidate for
\(\V[d'](x)\), and \(\V[d(x)](x)\) corresponds to the last suitable
candidate calculated under budget \(B\). Moreover, the expansion of the
tree is independent from the `leaf selection' step, and is fully
determined by the optimistic selection of leaves. We have the following
property for the returned leaf.

\begin{proposition}\label{finiteoptimum}
Given a budget $B\geq1$, Algorithm \ref{algo} terminates with output $S=\{\bm{u}_{\gamma,d(x)}^*(x),\V[d(x)](x)\}$  with horizon $d(x)\geq0$. $\hfill\Box$
\end{proposition}

\begin{IEEEproof}
Let  $x\in\mathbb{R}^n$  and $B\geq1$.   We  show  that  $S$ exactly  calculates  cost
$\V[d'](x)$ for  some $d'\in\mathbb{Z}_{>0}$.  The  optimal property  of output $S$ to Algorithm \ref{algo} is  fully determined  in the
particular iteration  in which  it is  updated.  Hence,  let $\mathcal{T}_i$  be the  tree to  be expanded  at iteration
$i\in\mathbb{Z}_{\geq 0}$, in which $S$ is  updated.  We show now that the selected leaf $S$ with  cost $J(S)$, where $J(S)$ is
the cost  associated to  leaf $S$,  attains the optimum  of horizon  $d':=\text{depth}(S)-1$, that  is $J(S)=\V[d'](x)$.
Since  $\V[d'](x)\leq J(S)$  by the  optimality of  $\V[d'](x)$, it  suffices to  prove $\V[d'](x)\geq  J(S)$. For  this
purpose, we proceed by contradiction,  and we assume that $\V[d'](x)<J(S)$. It follows from the fact that the input set $\mathcal{U}$ is finite that a
sequence that attains the optimum $\V[d'](x)$ exists, i.e.  there is a node $N\neq S$, descendent of root $x$ and possibly not
in  $\mathcal{T}_i$,  with cost  $J(N)=\V[d'](x)$.  Since  $\ell_u(x)\geq0$ for  any $x\in\mathbb{R}^n$  and
$u\in\{1,\ldots,M\}$, any  ancestor (parents, parents  of parents and  so on) of  $N$ will have  cost a lower  cost than
$J(N)$. Hence, let $L_i'$ be the ancestor of $N$ such that $L_i'\in\mathcal{L}(\mathcal{T}_i)$, thus $J(L_i')\leq J(N)$.
Then, we have $J(L_i')\leq J(N)=\V[d'](x)<J(S)$, that is $J(L_i')<J(S)$.  However, $S$ is the optimistically chosen leaf
$S=L_i$, and $J(S)=J(L_i)\leq J(L)$ for any leaf  $L\in\mathcal{L}(\mathcal{T}_i)$, hence for leaf $L_i'$,
it follows that $J(L_i')<J(S)\leq J(L_i')$: we have attained  a contradiction.  Therefore, $J(S)\leq\V[d'](x)$ and since $\V[d'](x)\leq J(S)$, we
conclude  $\V[d'](x)=J(S)$.  Thus,  at every  update of  $S$, a  new optimal  sequence is  found with  increased horizon
$d'\leftarrow d'+1$. Furthermore, note that at iteration $i=1$, the `leaf selection' step
is guaranteed to be entered and,  given budget  $B\geq1$, the outputs of Algorithm \ref{algo}, $d$ and $S$, are fully determined.
\end{IEEEproof}

The horizon \(d(x)\) in (\ref{eq:V}) depends on the given budget \(B\),
and will play a fundamental role in the near-optimality analysis
provided later in Section \ref{main-results}. The next proposition
provides a (conservative) relationship between budget \(B\) and a given
lower bound on \(d(x)\).

\begin{proposition}\label{dfiniteoptimum}
  Given $\bar d\geq 0$ and budget $B\geq\frac{M^{\bar d+1}-1}{M-1}$, Algorithm \ref{algo} returns $S=\{\bm{u}_{\gamma,d(x)}^*(x),\V[d(x)](x)\}$ with horizon $B-1\geq d(x)\geq\bar d$. $\hfill\Box$
\end{proposition}
\begin{IEEEproof}
The shallowest possible tree that can be explored with budget $B=\frac{M^{\bar d+1}-1}{M-1}$, under any circumstances, is the uniform, complete tree with depth $\bar d+1$, and one node of depth $\bar d+1$ expanded with children at depth $\bar d+2$. In this case $S$ is the node expanded at depth $\bar d+1$, hence $d(x) = \bar d$. Any other (e.g. optimistic) way of exploring the tree will lead to $d(x) \geq \bar d$. On the other hand, the deepest possible tree is the unbalanced tree, where the tree depth increase at every iteration, i.e. only one node is expanded per depth. Hence, given budget $B$, the maximum possible depth is $B$ and $d(x)=B-1$. Any other way of exploring the tree will
lead to $d(x)\leq B-1$. We conclude $B-1\geq d(x)\geq \bar d$.
\end{IEEEproof}

Proposition \ref{dfiniteoptimum} provides a relationship between a
minimum desired horizon \(\bar d\) in (\ref{eq:J}) and the required
budget to achieve it. This relationship is derived from the worst-case
exploration, which happens when OPmin is forced to uniformly explore the
possible choices of switches for a given horizon. Due to the optimistic
exploration, for given \(B=\frac{M^{\bar d+1}-1}{M-1}\), \(d(x)\) is
often much larger in practice than \(\bar d\). In the original OP study
\cite{hren2008}, this fraction is quantified by means of the branching
factor \(\kappa\), which we will investigate in future work.

\hypertarget{main-results}{%
\section{Main results}\label{main-results}}

In this section, we analyze the near-optimality properties of OPmin. We
also provide conditions under which system (\ref{eq:sys}), whose inputs
are generated in a receding-horizon fashion by OPmin, exhibits stability
properties. Assumptions are required for this purpose, which are now
stated.

\hypertarget{assumptions}{%
\subsection{Assumptions}\label{assumptions}}

We first assume that the optimization goal (\ref{eq:V}) is well-posed in
the following sense.

\begin{SA}
For any  $x \in  \mathbb{R}^n$, $\gamma  \in (0,1]$,  there exists  an infinite sequence of
admissible inputs $\bm{u}^*_{\gamma,\infty}(x)$, called \textit{optimal input sequence}, which minimizes (\ref{eq:Jinfty}), i.e. 
$\V[\infty](x)=\J[\infty](x,\ustar[\infty](x))$  is finite.  \hfill$\square$
\end{SA}

General conditions to ensure SA can be found in \cite{keerthi1985}.

We make the next general stabilizability and detectability assumptions
on system (\ref{eq:sys}) and stage cost \(\ell\) as in
\cite{grimm2005,romain2016,granzotto2019}, which are essential: (i) for
the construction of near-optimality bounds of the algorithm; (ii) to
ensure stability of the induced closed-loop system as demonstrated in
the sequel.

\begin{assumption}\label{Aone}
There exist $\overline{\alpha}_V,\alpha_W\in \mathcal{K}_\infty$, continuous functions  $W,\sigma: \mathbb{R}^n \to \mathbb{R}_{\geq 0}$, $\overline{\alpha}_W: \mathbb{R}_{\geq 0}\to
\mathbb{R}_{\geq 0}$ continuous, non-decreasing and zero at zero, such that the following conditions hold.

\begin{IEEEenumerate}[\settowidth{\labelwidth}{(ii)}]

 \item[(i)] For any $x \in \mathbb{R}^n$, $\gamma \in (0,1]$,
  \begin{equation} \V[\infty](x) \leq \overline{\alpha}_V(\sigma(x)).  \label{eq:avb} \end{equation}

 \item[(ii)] For any $x \in \mathbb{R}^n$, $u \in \mathcal{U}$,
   \begin{align*}
                W(x) &\leq \overline{\alpha}_W(\sigma(x))                   \addtocounter{equation}{1}\tag{\theequation}\label{eq:awb} \\
      W(f_u(x))-W(x) &\leq -\alpha_W(\sigma(x))+\ell_u(x).        \addtocounter{equation}{1}\tag{\theequation}\label{eq:aw}
   \end{align*}
 $\hfill\square$
\end{IEEEenumerate}
\end{assumption}

Function \(\sigma\) in Assumption \ref{Aone} serves as a measuring
function of the state and will be used to define stability, as in
\cite{romain2016,grimm2005,granzotto2019}. For instance, by defining
\(\sigma=|\cdot|\), one would be studying the stability of the origin,
and by taking \(\sigma=|\cdot|_{\cal A}\), one would study stability of
set \({\cal A}\subset \mathbb{R}^n\). Item (i) is related to the
stabilizability of system (\ref{eq:sys}) with respect to stage cost
\(\ell\). Indeed, it is shown in \cite[Lemma 1]{romain2016} that, if
stage cost \(\ell\) is uniformly globally exponentially controllable to
zero with respect to \(\sigma\) for system (\ref{eq:sys}), see
\cite[Definition 2]{grimm2005}, then Assumption \ref{Aone} is satisfied.
On the other hand, item (ii) of Assumption \ref{Aone} is a detectability
property of the stage cost \(\ell\) with respect to \(\sigma\). For
example, when \(\ell_u(x)=\sigma(x)\), one verifies item (ii) of
Assumption \ref{Aone} with \(W\equiv0\) and \(\alpha_W=\mathbb{I}\). For
a more general view on Assumption \(\ref{Aone}\), see the aforementioned
references. Note that we do not require \(\ell\) to take values in
\([0,1]\) contrary to \cite{hren2008}.

\hypertarget{relationship-between-vdxx-and-vinftyx}{%
\subsection{\texorpdfstring{Relationship between \(\V[d(x)](x)\) and
\(\V[\infty](x)\)}{Relationship between finite-horizon cost and infinite-horizon cost}}\label{relationship-between-vdxx-and-vinftyx}}

Algorithm \ref{algo} is able to calculate \(\V[d(x)](x)\) exactly for
any given \(x\in\mathbb{R}^n\), however it is not obvious how
\(\V[d(x)](x)\) relates to \(\V[\infty](x)\). Since \(\ell\) is not
constrained to take values in a given compact set, and we accept the
undiscounted case, the tools used in \cite{hren2008} to analyze
near-optimality are no longer applicable. We overcome this issue by
exploiting Assumption \ref{Aone}.

\begin{theorem} \label{Vestimates}
Suppose Assumption \ref{Aone} holds. For any $x  \in \mathbb{R}^n$, $\gamma\in(0,1]$ and $d(x)\in\mathbb{Z}_{>0}$, 
\begin{equation}\label{eq:Vestimates}
    \V[d(x)](x) \leq \V[\infty](x)           \leq \V[d(x)](x)+v_{\gamma,d(x)}(x),
\end{equation}
where  $v_{\gamma,d(x)}(x):=\gamma^{d(x)}\overline{\alpha}_V\circ\underline{\alpha}_Y^{-1}\circ\left(\frac{\mathbb{I}-{\alpha}_Y\circ\overline{\alpha}_Y^{-1}}{\gamma}\right)^{(d(x))}\circ\overline{\alpha}_Y(\sigma(x))$. Here,  ${\alpha}_Y=\underline{\alpha}_Y:=\alpha_W, \overline{\alpha}_Y:=\overline{\alpha}_V+\overline{\alpha}_W$ and $\alpha_W,\overline{\alpha}_W,\overline{\alpha}_V$ come from Assumption \ref{Aone}. $\hfill\square$
\end{theorem}

Theorem \ref{Vestimates} follows from \cite[Theorem 3]{granzotto2019},
therefore the proof is omitted. The lower-bound in (\ref{eq:Vestimates})
trivially holds from the optimality of \(\V[d(x)](x)\) as
\(d(x)<\infty\), see Proposition \ref{dfiniteoptimum}. The upper-bound,
on the other hand, implies that the infinite-horizon cost is at most
\(v_{\gamma,d(x)}(x)\) away from the finite-horizon \(\V[d(x)](x)\). The
error term \(v_{\gamma,d(x)}(x)\) has three desirable properties
compared to the term given in \cite{hren2008}, which we recall is
\(\frac{\gamma^{d(x)}}{1-\gamma}\). First, when \(\sigma(x)\) is small,
so is \(v_{\gamma,d(x)}(x)\). Second, \(v_{\gamma,d(x)}(x)\) is finite
for \(\gamma=1\), while in \cite{hren2008},
\(\frac{\gamma^{d(x)}}{1-\gamma}\to\infty\) in this case. Third,
\(v_{\gamma,d(x)}(x)\to 0\) when \(d(x)\to\infty\) for \(\gamma\)
sufficiently close to 1, as seen in \cite[Lemma 3]{granzotto2019}, which
is true for OP only when \(\gamma<1\). Thus, in contrast to OP, by
exploiting stabilizability and detectability properties, we have
obtained an error bound that forfeits the assumption \(\ell\in[0,1]\),
accepts the undiscounted case \(\gamma=1\), and is decreasing in
\(d(x)\), even when \(\gamma=1\).

\begin{remark}
Lemma 3 in \cite{granzotto2019} states that, when $\gamma\in(1-\frac{{\alpha}_Y(\Delta)}{\overline{\alpha}_Y(\Delta)},1]$
and $\sigma(x)\leq\Delta$ for any given $\Delta\geq 0$,  the error $v_{\gamma,d(x)}(x)\to 0$ when $d(x)\to\infty$. In fact, it is shown in the proof of \cite[Lemma 3]{granzotto2019} that for $\gamma\in(1-\frac{{\alpha}_Y(\Delta)}{\overline{\alpha}_Y(\Delta)},1]$, $v_{\gamma,d(x)}\leq\gamma^{d(x)}\overline{\alpha}_V\circ\underline{\alpha}_Y^{-1}(\overline{\alpha}_Y(\Delta))$, which has the same decrease rate
in $d(x)$ as the original OP error bound $\gamma^{d(x)}\frac{1}{1-\gamma}$. However
for $\gamma=1$ it follows that $v_{\gamma,d(x)}$ decreases as a function of $(\mathbb{I}-{\alpha}_Y\circ\overline{\alpha}_Y^{-1})^{(d(x))}$, which is not the case for  OP  under the assumptions of \cite{hren2008}. We can then control the size of $v_{\gamma,d(x)}$ as small as desired by chosing a suitable $d(x)$, via a suitable budget,   and $\gamma$.  $\hfill\Box$
\end{remark}

It is unclear if increasing or decreasing \(\gamma\) for a fixed
\(d(x)\) will increase or decrease error bound \(v_{\gamma,d(x)}(x)\).
Indeed, this is due to competing terms \(\gamma^{d(x)}\), which
increases as \(\gamma\) increases for fixed \(d(x)\), and
\(\left(\frac{\mathbb{I}-{\alpha}_Y\circ\overline{\alpha}_Y^{-1}}{\gamma}\right)^{(d(x))}\),
which decreases as \(\gamma\) increases with fixed \(d(x)\). When
Assumption \ref{Aone} is satisfied with class \(\mathcal{K}_\infty\)
functions of a particular form, we show that the error term is uniform
in \(\gamma\), thus clarifying this issue.

\begin{corollary}\label{Vestimatesunif}
Suppose  that Assumption \ref{Aone} is satisfied and there exist $\bar a_W \geq 0$, $a_W,\bar a_V > 0$ such that $\overline{\alpha}_V(s)
\leq \bar a_V \cdot s$, $\overline{\alpha}_W(s)  \leq \bar a_W \cdot s$, $\alpha_W(s) \geq a_W \cdot s$  for any $s \geq
0$. Let $x \in  \mathbb{R}^n$, $d(x)\in\mathbb{Z}_{>0}$,  any $\gamma\in(0,1]$,
            \begin{equation}\V(x) \leq  \V[\infty](x)\leq \V(x)+\hat v_{d(x)}(x)
                                          \label{eq:VestimatesunifA}
            \end{equation}
            where  $\hat v_{d(x)}(x):= \displaystyle\frac{\bar a_V(\bar a_V+\bar a_W)}{a_W}\left(1-\frac{a_W}{\bar a_V+\bar a_W}\right)^{d(x)}\!\!\!\!\sigma(x)$.

$\hfill\square$
\end{corollary}

The proof of Corollary \ref{Vestimatesunif} is a direct substitution of
the linear terms in Theorem \ref{Vestimates}, and is therefore omitted.
Compared to Theorem \ref{Vestimates}, Corollary \ref{Vestimatesunif}
provides an error bound \(\hat v_{d(x)}(x)\) uniform in \(\gamma\),
linear in \(\sigma(x)\), which also decreases exponentially to 0 for any
\(\gamma\in(0,1]\) when \(d(x)\to\infty\), as a
function\footnote{Indeed,  $\left(1-\frac{a_W}{\bar a_V+\bar a_W}\right)\in[0,1)$ holds as $a_W s\leq \alpha_W(s)\leq\overline{\alpha}_V(s)+\overline{\alpha}_W(s)\leq(\bar a_V+\bar a_W)s$ for any $s\geq0$, in view of Corollary \ref{Vestimatesunif} and Proposition \ref{YLyapunovProp}.}
of \(\left(1-\frac{a_W}{\bar a_V+\bar a_W}\right)^{d(x)}\).

\hypertarget{stability}{%
\subsection{Stability}\label{stability}}

We now consider the scenario where (\ref{eq:sys}) is controlled in a
receding horizon fashion by OPmin as defined by Algorithm \ref{algo}.
That is, at each time instant \(k\in\mathbb{Z}_{\geq 0}\), the first
element of the optimal sequence \(\bm{u}^*_{\gamma,d(x_k)}(x_k)\), is
calculated by OPmin, and then applied to system (\ref{eq:sys}). This
leads to closed-loop system
\begin{equation}   x_{k+1}  \in  f_{\mathcal{U}^*_{\gamma,d(x_k)}(x_k)}(x_k)  =:
F^*_{\gamma,d(x_k)}(x_k), \label{eq:autosys}\end{equation} where
\(f_{\mathcal{U}^*_{\gamma,d(x)}(x)}(x)\) is the set
\(\{f_u(x) : u \in \mathcal{U}^*_{\gamma,d(x)}(x)\}\), and
\(\mathcal{U}_{\gamma,d(x)}^*(x):= \big\{ u_0 : \exists u_1,\ldots,u_{d(x)} \in \mathcal{U}\text{ such that }  \V[d(x)](x)=\J[d(x)](x,[u_0,\ldots,u_{d(x)}])\big\}\)
is the set of the first input of \(d(x)\)-horizon optimal input
sequences at \(x\). We denote by \(\phi(k,x)\), with some abuse of
notation, a solution to (\ref{eq:autosys}) at time
\(k\in\mathbb{Z}_{\geq 0}\) with initial condition \(x\in\mathbb{R}^n\).

The next theorem provides stability guarantees for system
(\ref{eq:autosys}).

\begin{theorem}\label{algostab}
Consider system (\ref{eq:autosys}) and  suppose Assumption \ref{Aone} holds. There exists $\beta \in  \mathcal{KL}$ such that for
any  $\delta,\Delta>0$,  there  exist $\gamma^*  \in  (0,1)$  and  $\bar d\in\mathbb{Z}_{>0}$  such  that for  any  $\gamma  \in
(\gamma^*,1]$, any  budget  $B\geq\frac{M^{\bar d+2}-1}{M-1}$, any $x  \in  \{z  \in\mathbb{R}^n  \, :  \,  \sigma(z)  \leq  \Delta \}$,  any  solution
$\phi(\cdot,x)$ to system (\ref{eq:autosys}) satisfies, for all $k\in
\mathbb{Z}_{\geq 0}$

\begin{equation}
  \sigma(\phi(k,x)) \leq \max \{\beta(\sigma(x),k),\delta \}. \label{eq:Ystab}
\end{equation}~$\hfill\square$
\end{theorem}

The proof of Theorem \ref{algostab} is given in the appendix. Theorem
\ref{algostab} provides a semiglobal practical stability property for
set \(\{z : \sigma(z) = 0 \}\). This implies that solutions to
(\ref{eq:autosys}), with initial state \(x\) such that
\(\sigma(x)\leq\Delta\), where \(\Delta\) is any given (arbitrarily
large) strictly positive constant, will converge to the set
\(\{z : \sigma(z) \leq \delta\}\), where \(\delta\) is any given
(arbitrarily small) strictly positive constant, by taking \(\gamma\)
sufficiently close to \(1\) and a budget sufficiently large. Note that
we take budget \(B\geq\frac{M^{\bar d+2}-1}{M-1}\) to guarantee
\(d(x)>\bar d\), by Proposition \ref{dfiniteoptimum}. While OP was used
in various control problems
\cite{BUSONIU2017297,busoniunetwork,7963352}, such stability properties
were never been proved before.

By strengthening Assumption \ref{Aone}, we can prove a global
exponential stability property. Since it follows the same arguments as
in \cite[Corollary 2]{granzotto2019} and the modifications given in the
appendix, the proof is omitted.

\begin{corollary}\label{Yges}
Suppose that the conditions of Corollary  \ref{Vestimatesunif} holds.  Let  $\gamma^*, \bar d$  be such  that
\begin{equation}
 1-\gamma^*+\frac{\bar  a_V}{a_W}\left(1-\frac{a_W}{\bar  a_V  +\bar  a_W}\right)^{\bar d}  <
   \frac{a_W}{\bar a_V+\bar  a_W}. \label{eq:gdcondges}
\end{equation}
Then, there exist  $K,\lambda>0$, such that for  any $\gamma\in(\gamma^*,1]$, any budget $B\geq\frac{M^{\bar d+2}-1}{M-1}$, for any $x  \in \mathbb{R}^n$, the solution  $\phi(\cdot,x)$ to system (\ref{eq:autosys}),  satisfies $\sigma(\phi(k,x)) \leq
K\sigma(x)e^{-\lambda k}$ for all $k \in \mathbb{Z}_{\geq 0}$. $\hfill\square$
\end{corollary}

Corollary \ref{Yges} ensures a uniform global exponential stability
property of \(\{x:\sigma(x)=0\}\) for (\ref{eq:autosys}). Inequality
(\ref{eq:gdcondges}) is always feasible for \(\gamma^*\) sufficiently
close to 1 and \(\bar d\) sufficiently large. Indeed, we either first
fix \(\gamma^*\in(\bar\gamma,1]\) with
\(\bar\gamma=1-\tfrac{a_W}{\bar a_V+\bar a_W}\) and then select
\(\bar d\) and the associated budget \(B\) such that
(\ref{eq:gdcondges}) holds, or we first fix budget \(B\) with associated
\(\bar d> \tilde d\) with
\(\tilde d=\floor{\tfrac{\ln(\bar a_V (\bar a_V +\bar a_W)/a_W^2)}{-\ln (1-\frac{a_W}{\bar a_V+\bar a_W})}}\)
and select \(\gamma^*\) such that (\ref{eq:gdcondges}) holds. The
resulting pair \((\gamma^*,\bar d)\) and budget \(B\) are suitable
candidates for (\ref{eq:gdcondges}) by construction. This is consistent
with results for finite-horizon discounted costs \cite{granzotto2019},
where both \(\gamma\) has to be sufficiently close to 1 and \(\bar d\)
has to be sufficiently large, and results for finite-horizon
undiscounted cost \cite{grimm2005}, where \(d\) has to be taken large.

\begin{remark}
It is possible to relax Corollary 2 conditions and derive semiglobal asymptotic results, similarly to \cite[Corollary 1]{granzotto2019}. $\hfill\Box$
\end{remark}

\hypertarget{near-optimality-guarantees}{%
\subsection{Near-optimality
guarantees}\label{near-optimality-guarantees}}

In Theorem \ref{Vestimates}, we have provided near-optimality guarantees
of finite-horizon cost \(\V[d(x)]\) with respect to the infinite-horizon
cost \(\V[\infty](x)\). This is an important feature of OPmin, but this
does not directly provide us with information on the actual value of the
cost function (\ref{eq:Jinfty}) along solutions to (\ref{eq:autosys}).
Indeed, we do not implement the whole sequence
\(\bm{u}^*_{\gamma,d(x)}(x)\) given by OPmin at \(x\) in
(\ref{eq:autosys}), instead we proceed in a receding horizon fashion.
The relevant cost function to analyze is thus the running cost
\cite{gruneperformance} defined as \begin{equation}
        \mathcal{V}_{\gamma,\bar d}^{\text{run}}(x) := \sum_{k=0}^\infty \gamma^k\ell_{\mathcal{U}^*_{\gamma,d(\phi(k,x))}(\phi(k,x))}(\phi(k,x)),
 \label{eq:Vrun}
\end{equation} where \(d(\phi(k,x))>\bar d\) for all
\(k\in\mathbb{Z}_{\geq 0}\), and \(\bar d\) is a lower bound on the
desired horizon at each step, which we can enforce by taking a
sufficiently large budget according to Proposition \ref{dfiniteoptimum}.
It has to be noted that \(\mathcal{V}^\text{run}_{\gamma,\bar d}(x)\) is
a set, since solutions of (\ref{eq:autosys}) are not necessarily unique.
Each element
\(V_{\gamma,\bar d}^{\text{run}}(x) \in \mathcal{V}_{\gamma,\bar{d}}^{\text{run}}(x)\)
corresponds then to the cost of a solution of (\ref{eq:autosys}).
Clearly, \(\V[\bar d]^{\text{run}}(x)\) is not necessarily finite, as
the stage costs may not decrease to 0 in view of Theorem \ref{algostab}.
Indeed, only practical convergence is ensured in Theorem \ref{algostab}
in general. As a result, the corresponding running cost may not be
finite. We therefore restrict our attention to the case where Corollary
\ref{Yges} holds, in the next theorem.

\begin{theorem}\label{Vrunestimates}
Consider system (\ref{eq:autosys}) and assume that Corollary \ref{Yges} holds with tuple $(K,\lambda,\gamma^*,\bar d)$. For any $\gamma\in(\gamma^*,1]$, budget $B\geq\frac{M^{\bar d+2}-1}{M-1}$, $x  \in \mathbb{R}^n$, and $V_{\gamma,\bar d}^{\text{run}}(x) \in \mathcal{V}_{\gamma,\bar d}^{\text{run}}(x)$, 
\begin{equation}
\V[\infty](x)\leq\V[\bar d]^{\text{run}}(x) \leq     \V[\infty](x) + w_{\gamma,\bar d} \cdot\sigma(x),
\label{eq:Vrunestimates}
\end{equation}
where $
w_{\gamma, \bar d}:=\displaystyle\left(1-\frac{a_W}{\bar a_V+\bar a_W}\right)^{\bar d}\frac{K\bar a_V(\bar a_V+\bar a_W)\gamma}{a_W(e^{\lambda}-\gamma)}.
$
$\hfill\square$
\end{theorem}

\begin{IEEEproof}
Let  $x   \in  \mathbb{R}^n$,   $\gamma  \in   (\gamma^*,1]$, budget $B\geq\frac{M^{\bar d+2}-1}{M-1}$, and  $\phi(k+1,x)\in F^*_{\gamma,d(\phi(k,x))}(\phi(k,x))$ for any $k\in\mathbb{Z}_{\geq 0}$ where $\phi$ is a solution to (\ref{eq:autosys}) initialized at $x$. For the sake of convenience, we denote $\ell(x,u):=\ell_u(x)$ for any $x\in\mathbb{R}^n$ and $u\in\mathcal{U}$. Consider 
\begin{equation}
  V_{\gamma, \bar d}^\text{run}(x):=\sum_{k=0}^{\infty}\gamma^k\ell(\phi(k,x),u^r_k) \label{eq:Vrunphi},
\end{equation} 
where $u^r_k\in\mathcal{U}^*_{\gamma,d(\phi(k,x))}(\phi(k,x))$ such that $\phi(k+1,x)=f_{u^r_k}(\phi(k,x))$. Note that indeed  $V_{\gamma, \bar d}^\text{run}(x)\in\mathcal{V}_{\gamma,\bar d}^{\text{run}}(x)$. The inequality $\V[\infty](x) \leq V^{\text{run}}_{\gamma, \bar d}(x)$  follows from optimality of $\V[\infty](x)$. Since $u^r_k\in\mathcal{U}^*_{\gamma,d(\phi(k,x))}(\phi(k,x))$, we derive from Bellman optimality principle that 
\begin{equation}
\begin{split}
\ell(\phi(k,x),u^r_k)&=\V[d(\phi(k,x))](\phi(k,x))\\
                     &\quad-\gamma\V[d(\phi(k,x))-1](\phi(k+1,x)).
\end{split} \label{eq:bellmanA}
\end{equation} Since, for any $z\in\mathbb{R}^n$, $d(z)>\bar d$ holds
for  budget $B$  according to Proposition \ref{dfiniteoptimum}, $d(\phi(k,x))-1\geq\bar d$ and $\V[d(\phi(k,x))-1](\phi(k+1,x))\geq\V[\bar d](\phi(k+1,x))$ for all $k\in\mathbb{Z}_{\geq 0}$. As a result, in view of (\ref{eq:bellmanA}), 
$\ell(\phi(0,x),u^r_0)\leq\V[d(\phi(0,x))](\phi(0,x))-\gamma\V[\bar d](\phi(1,x))$. On the other hand, since $\V[d(\phi(k,x))](\phi(k,x))\leq\V[\infty](\phi(k,x))$,
$\ell(\phi(k,x),u^r_k)\leq\V[\infty](\phi(k,x))-\gamma\V[\bar d](\phi(k+1,x))$.  Thus
\begin{equation}
\begin{split}
  V_{\gamma, \bar d}^\text{run}(x)&\leq\V[d(\phi(0,x))](\phi(0,x)) -\gamma\V[\bar d](\phi(1,x))\\
                               &\quad +\gamma\ (\V[\infty](\phi(1,x))-\gamma\V[\bar d](\phi(2,x)))\\
                               &\quad +\gamma^2(\V[\infty](\phi(2,x))-\gamma\V[\bar d](\phi(3,x)))+\ldots\\
                               &=\V[d(\phi(0,x))](\phi(0,x)) \\
                               &\quad + \sum_{k=1}^{\infty}\gamma^k(\V[\infty](\phi(k,x))-\V[\bar d](\phi(k,x))).
\end{split}\label{eq:Vrundiffphi}
\end{equation}
According to  Corollary \ref{Vestimatesunif},  $\V[\infty](\phi(k,x))-\V[\bar d](\phi(k,x))\leq \hat v_{\bar d}(\phi(k,x))$, with $\hat v_{\bar d}(z)=\frac{\bar a_V}{a_W}\left(1-\frac{a_W}{\bar a_V+\bar a_W}\right)^{\bar d}(\bar a_V+\bar a_W)\sigma(z)$ for any $z\in\mathbb{R}^n$. 
Hence, by direct substitution in (\ref{eq:Vrundiffphi}), $\V[\bar d]^{\text{run}}(x) \leq     \V[d(\phi(0,x))](\phi(0,x)) \quad + \tfrac{\bar a_V(\bar a_V+\bar a_W)}{a_W}\left(1-\tfrac{a_W}{\bar a_V+\bar a_W}\right)^{\bar d}\sum_{k=1}^\infty \gamma^k \sigma(\phi(k,x))$.
Recalling  $\phi(0,x)=x$ and that $\sigma(\phi(k,x)) \leq K\sigma(x)e^{-\lambda k}$ holds from Corollary \ref{Yges}, we obtain
\begin{equation}
\begin{split}
\V[\bar d]^{\text{run}}(x) &\leq     \V[d(x)](x) \\
                           &\quad + \sigma(x)\tfrac{K\bar a_V(\bar a_V+\bar a_W)}{a_W}\left(1-\tfrac{a_W}{\bar a_V+\bar a_W}\right)^{\bar d}\sum_{k=1}^\infty \gamma^k e^{-\lambda k},\\
                           &\leq    \V[\infty](x) \\
                           &\quad + \sigma(x)\tfrac{K\bar a_V(\bar a_V+\bar a_W)}{a_W}\left(1-\tfrac{a_W}{\bar a_V+\bar a_W}\right)^{\bar d}\frac{\gamma}{e^{\lambda}-\gamma}.
\end{split}
\label{eq:Vrundiffsum}
\end{equation}
Since (\ref{eq:Vrundiffsum}) holds
for an arbitrary solution of (\ref{eq:autosys}), $\phi(k+1,x)=f_{u^r_k}(\phi(k,x))$ for any $k\in\mathbb{Z}_{\geq 0}$,  (\ref{eq:Vrundiffsum}) holds for any $V_{\gamma, \bar d}^\text{run}(x)\in\mathcal{V}_{\gamma,\bar d}^{\text{run}}(x)$.
\end{IEEEproof}

Similarly to Theorem \ref{Vestimates}, the inequality
\(\V[\infty](x) \leq V^{\text{run}}_{\gamma, \bar d}(x)\) of Theorem
\ref{Vrunestimates} directly follows from the optimality of
\(\V[\infty](x)\). On the other hand, the inequality
\(\V[\bar d]^{\text{run}}(x) \leq \V[\infty](x) + w_{\gamma,\bar d}(x)\)
provides a relationship between the running cost
\(V_{\gamma,\bar d}^\text{run}(x)\) and the infinite-horizon cost at
state \(x\), \(\V[\infty](x)\). The term \(w_{\gamma,\bar d}\) can be
explicitly calculated, see Corollary \ref{Vestimatesunif} and
\cite[proof of Corollary 3]{granzotto2019} for the expressions of \(K\)
and \(\lambda\). The latter inequality in (\ref{eq:Vrunestimates})
confirms the intuition coming from Theorem \ref{Vestimates} that a large
computational budget \(B\) leads to tight near-optimality guarantees.
That is, when \(\bar d\to\infty\), or equivalently \(B\to\infty\)
according to Proposition \ref{dfiniteoptimum}, \(w_{\gamma,\bar d}\to0\)
and \(V^{\text{run}}_{\gamma,d}(x)\to\V[\infty](x)\), provided that
\(\gamma\) and budget \(B\) have been chosen as to stabilize system
(\ref{eq:autosys}). In contrast with Theorem \ref{Vestimates}, stability
of system (\ref{eq:autosys}) plays a role in Theorem
\ref{Vrunestimates}. Indeed, the term
\(K\frac{\gamma}{e^{\lambda}-\gamma}\) in (\ref{eq:Vrunestimates}) shows
that the larger the exponential decay \(\lambda\) is, the smaller the
error term \(w_{\gamma,\bar d}\) will be. The running cost for the
original OP was considered in \cite{busoniunetwork}, and it was found to
perform at worst like the finite sequence,
i.e.~\(\V[\bar d]^{\text{run}}(x) \leq \V[\infty](x) + \frac{\gamma^{\bar d}}{1-\gamma}\).
Compared to the bound derived for OP, the bound in Theorem
\ref{Vrunestimates} has similar benefits as Theorem \ref{Vestimates},
namely we are not limited to \(\ell\in[0,1]\), it does not explode to
\(\infty\) when \(\gamma\) tends to 1, and when \(\sigma(x)\) is small
follows \(w_{\gamma,\bar d}\cdot\sigma(x)\) small. Moreover, the
mismatch decays exponentially in \(\bar d\), independent on \(\gamma\).

\begin{remark}
Inequality (\ref{eq:Vrunestimates}) can be written as a relationship of the finite-horizon costs in view of Theorem \ref{Vrunestimates}. In particular, we have
$
\V[d(x)](x)\leq\V[\infty](x)\leq\V[\bar d]^{\text{run}}(x) \leq     \V[d(x)](x) + w_{\gamma,\bar d}\cdot\sigma(x)
$, for any $x\in\mathbb{R}^n$.
Hence, $\V[d(x)](x)$ can be used to upper and lower bound $\V[\bar d]^{\text{run}}(x)$ from the first call of OPmin at initial state $x$. $\hfill\Box$
\end{remark}

Other authors have considered the running cost of finite-horizon
controllers applied in receding horizon fashion, like
\cite{gruneperformance} in the context of model predictive control. In
particular, \cite{gruneperformance} derives relative performance of the
running cost to the infinite-horizon optimal cost. From Theorem
\ref{Vrunestimates}, we derive a similar result.

\begin{corollary}\label{Vrunrelative}
Suppose Theorem \ref{Vrunestimates} holds for system (\ref{eq:autosys}) with tuple $(K,\lambda,\gamma^*,\bar d)$. Then, for any $x\in\mathbb{R}^n$, $\gamma\in(\gamma^*,1]$, budget $B\geq\frac{M^{\bar d+2}-1}{M-1}$, and $V_{\gamma,\bar d}^{\text{run}}(x) \in \mathcal{V}_{\gamma,\bar d}^{\text{run}}(x)$ such that $\V[\infty](x)>0$, 
\begin{equation}
\frac{\V[\bar d]^{\text{run}}(x)-\V[\infty](x)}{\V[\infty](x)+W(x)}\leq \frac{w_{\gamma,\bar d}}{a_W}.
\end{equation}$\hfill\Box$
\end{corollary}
\begin{IEEEproof}[Sketch of proof]
Let  $x   \in  \mathbb{R}^n$,   $\gamma  \in   (\gamma^*,1]$, budget $B\geq\frac{M^{\bar d+2}-1}{M-1}$. The proof follows by substitution of $a_W\sigma(x)\leq \V[\infty](x)+W(x)$ from item (i) of Proposition \ref{YLyapunovProp} in the appendix and $\V[\bar d](x)\leq\V[\infty](x)$ in
 $\V[\bar d]^{\text{run}}(x) \leq     \V[\infty](x) + w_{\gamma,\bar d}\cdot\sigma(x)$ from Theorem \ref{Vrunestimates}.
\end{IEEEproof}

Corollary \ref{Vrunrelative} provides the relative performance of
\(\V[\bar d]^{\text{run}}(x)+W(x)\) and \(\V[\infty](x)+W(x)\). When
Assumption \ref{Aone} holds with \(W\equiv0\), Corollary
\ref{Vrunrelative} provides a relative relationship between any running
cost \(\V[\bar d]^{\text{run}}\) and the infinite-horizon cost
\(\V[\infty]\), as done in \cite{gruneperformance}. Interestingly, the
obtained bound for relative performance is uniform in \(x\). It does
however conserve the desired properties of Theorem \ref{Vrunestimates},
in particular the exponential decay in \(\bar d\) in view of the
expression of \(w_{\gamma,\bar d}\) in Theorem \ref{Vrunestimates}.

Corollary \ref{Vrunrelative} follows from the stabilizability and
detectability properties in Assumption 1, when the closed-loop system
satisfies global exponential stability, while \cite{gruneperformance}
results derive from relaxed dynamic programming, which relies on
parameters of a modified Bellman optimality equation. Although a direct
comparison between results is thus difficult, we note that our results
are exponential in minimal horizon \(\bar d\) with rate
\(\left(1-\frac{a_W}{\bar a_V+\bar a_W}\right)^{\bar d}\), while in
\cite{gruneperformance} the obtained bound is of order
\(\frac{1}{(\eta+1)^N-1}\) for a
parameter\footnote{We use $\eta$ instead of $\gamma$, as used in \cite{gruneperformance}, to not be confused with our discount factor $\gamma$.}
\(\eta\) derived from the relaxed dynamic programming property and
constant horizon \(N+1\). Hence, we expect that our bound may provide
similar relative performance with a shorter horizon. We note that, while
in \cite{gruneperformance} stability of the system is not assumed as in
Corollary \ref{Vrunrelative}, a lower-bound for horizon \(N+1\) such
that the results of \cite{gruneperformance} holds is expected, as in
Corollary \ref{Vrunrelative}.

\hypertarget{example}{%
\section{Example}\label{example}}

We consider the cubic integrator from \cite[Example 1]{grimm2005},
i.e.~\(x_1^+=x_1+u\), \(x_2^+=x_2+u^3\) where
\((x_1,x_2):=x\in\mathbb{R}^2\) and \(u\in\mathbb{R}\). It was verified
in \cite{grimm2005} that an open-loop sequence of inputs drives the
system to \(x=0\) in a finite number of steps. This open-loop sequence
can be expressed as three feedback gains \(K_1(x)=-x_1\),
\(K_2(x)=x_2^{\frac{1}{3}}\) and
\(K_3(x)=\left(-\frac{1}{2}+\sqrt{\frac{7}{12}}\right)x_2^{\frac{1}{3}}\),
which are successively applied. We propose here to switch between these
gains to minimize cost (\ref{eq:Jinfty}), with
\(\ell_u(x)=|x_1|^3+|x_2|+|K_u(x)|^3\) for any \(x\in\mathbb{R}^2\) and
\(u\in\{1,2,3\}\). Note that we cannot design a local LQR controller for
this system, due the lack of stabilizability of the linearized model at
the origin. We therefore consider the switched system
\(x_1^+=x_1+K_u(x)\), \(x_2^+=x_2+(K_u(x))^3\), for \(u\in\{1,2,3\}\).
We apply OPmin to illustrate the near-optimality and stability
guarantees of Section \ref{main-results}. To do this, note that SA
applies for the same reasons as in \cite{grimm2005}. By taking
\(\sigma(x)=|x_1|^3+|x_2|\) for any \(x\in\mathbb{R}^2\), Assumption
\ref{Aone} holds \(\alpha_W=\mathbb{I}\), \(W=\alpha_W=0\) and
\(\overline{\alpha}_V=14\mathbb{I}\), as in \cite{grimm2005}. We verify
Corollary \ref{Yges} conditions with \(a_W=1\), \(\bar a_V=14\) and
\(\bar a_W=0\), and conclude that for \(\gamma=1\), any budget
\(B\geq \frac{3^{73}-1}{2}\) ensures global exponential stability.
Consequently, Theorem \ref{Vrunestimates} also holds. The lower-bound on
\(B\) is conservative, as the horizon \(\bar d\) itself in Corollary
\ref{Yges} is subject to some conservatism, and that OPmin will ensure
large horizons for smaller budgets in general. We have thus fixed the
budget to \(B=3000\) for initial condition \(x=[-1,1.5]^\top\). Figure 1
shows the evolution of the state, and we see that both \(x_1\) and
\(x_2\) converge to zero, as ensured by Corollary \ref{Yges}. We then
consider several initial conditions and we study the impact of the
budget on the actual running cost estimated by running simulations over
200 steps. We see in Table I that the estimated running cost becomes
smaller when increasing the budget, which is consistent with Theorem
\ref{Vrunestimates}. In other words, the larger the budget, the better
the running cost performance.

\vspace{-1em}
\begin{figure}[H]
\centering
\includegraphics[width=0.40\textwidth]{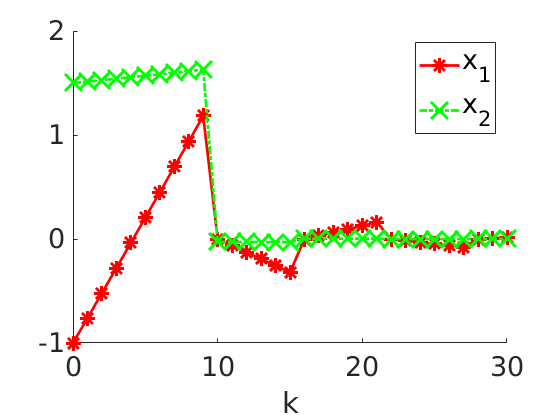}
\vspace{-0.3em}
\caption{State evolution for $B=3000$ and $x=[-1,1.5]^\top$.}
\end{figure}
\vspace{-2em}
\begin{table}[H]
\begin{center}
\begin{tabular}{
 r 
 r 
 |c c c 
}
\bottomrule
&       & \multicolumn{1}{c}{\vphantom{\Big(}}&{Budget}&\\[0.3em]
&   \    & \multicolumn{1}{c}{30} & \multicolumn{1}{c}{300} & \multicolumn{1}{c}{3000}\\
\hline
                                 &$\vphantom{\Big(}[\hphantom{-}10,\hphantom{.}\hphantom{-}15]^{\top}$            &  199015 & 13757 &  12609\\
\multirow{2}{*}{\begin{minipage}{2em}Initial States\end{minipage}}       &$\vphantom{\Big(}[-1\hphantom{0},\hphantom{-}1.5]^{\top}$    &  314  &  28  &  22\\
                                 &$\vphantom{\Big(}[-15,\, -10]^{\top}$   &  128184477 & 46875 &  42952\\
                                 &$\vphantom{\Big(}[\hphantom{-}10,\, -15]^{\top}$           &  14180  & 2802  &  2615\\[-0.25em]
\bottomrule
\end{tabular}
\caption{Estimated running cost for various budgets and initial conditions.}
\end{center}
\end{table}

\vspace{-3em}

\hypertarget{conclusion}{%
\section{Conclusion}\label{conclusion}}

We have modified the optimistic planning algorithm in \cite{hren2008} to
be applicable for the near-optimal, stable control of nonlinear switched
discrete-time systems. We relied for this purpose on general
stabilizability and detectability assumptions, originally stated in the
model predictive control literature \cite{grimm2005}. We have then
analyzed the algorithm near-optimality guarantees, which has major
features over the bound in \cite{hren2008} as discussed in Section
\ref{main-results}. We have also shown that a system controlled in a
receding-horizon fashion by OPmin satisfies stability properties. We
have finally analyzed the mismatch between the optimal value function
and the obtained running cost, and the same benefit as for the
near-optimality guarantees were observed.

\renewcommand\thesection{Appendix}

\hypertarget{proof-of-theorem}{%
\section{\texorpdfstring{Proof of Theorem \ref{algostab}
\label{proof}}{Proof of Theorem  }}\label{proof-of-theorem}}

The proof of Theorem \ref{algostab} follows the same steps as the proof
of \cite[Theorem 2]{granzotto2019}. The difference is that the horizon
in cost (\ref{eq:V}) is not fixed as in \cite{granzotto2019}, but
depends on the state. Nevertheless, as noted in
\cite[Remark 3]{granzotto2019}, the results can be modified to hold for
varying horizons, provided that considered horizons are lower-bounded by
a sufficiently large constant, \(\bar d\) in our case. This is what we
explicitly show in the following. We first state the next Lyapunov
properties.

\begin{proposition}\label{YLyapunovProp}
Suppose Assumption  \ref{Aone} holds. For  any $\gamma \in  (0,1]$ and
$\bar d\in\mathbb{Z}_{>0}$,  there exists $\Y[d(\cdot)] : \mathbb{R}^n \to \mathbb{R}_{\geq 0}$ such that the following holds.
\begin{IEEEenumerate}[\settowidth{\labelwidth}{(ii)}]
   \item[(i)] For any $x\in\mathbb{R}^n$, $ \underline{\alpha}_Y(\sigma(x)) \leq \Y[d(x)](x)  \leq \overline{\alpha}_Y(\sigma(x)),$
where $\underline{\alpha}_Y,\overline{\alpha}_Y$ come from Theorem \ref{Vestimates}.
   \item[(ii)] For any $x\in\mathbb{R}^n$ with $d(x)\geq\bar d$, $v \in F^*_{\gamma,d(x)}(x)$, $ \Y[d(v)](v)-\Y[d(x)](x) \leq \frac{1}{\gamma}\Big(-{\alpha}_Y(\sigma(x)) + \Upsilon(\Y[d(x)](x),\gamma,\bar d)\Big)$
     where $\Upsilon(s,\gamma,\bar d):=(1-\gamma)s+\gamma^{\bar d}\overline{\alpha}_V\circ\underline{\alpha}_Y^{-1}\circ\left(\frac{\mathbb{I}-{\alpha}_Y\circ\overline{\alpha}_Y^{-1}}{\gamma}\right)^{(\bar d)}(s)$, where $\overline{\alpha}_V$ comes from Assumption \ref{Aone} and $\alpha_Y$ is defined in Theorem \ref{Vestimates}, and for any $s\geq0$, $\Upsilon(s,\gamma,\bar d)\to0$ when $\gamma\to1$ and $\bar d\to\infty$. $\hfill\square$
\end{IEEEenumerate}
\end{proposition}

\begin{IEEEproof}
The proof works by following the steps in \cite[proof of Theorem 1]{granzotto2019}, only one step needs to be carefully modified. To show the mechanism, we concentrate on the case where $\gamma<1$; the same reasoning applies when $\gamma=1$.
Let  $\gamma  \in  (0,1)$,  $\bar d  \in  \mathbb{Z}_{>0}$,  $x\in\mathbb{R}^n$  and  $v  \in  F^*  _  {\gamma,d(x)}(x)$, with $d(x)\geq \bar d$ and $d(v)\geq \bar d$. 
There  exists
$[\us[0],\us[1],\ldots,\us[d(x)]]=\bm{u}^*_{\gamma,d(x)}(x)$ such  that $v=f_{\us[0]}(x)$ and  $\bm{u}^*_{\gamma,d(x)}(x)$ is  an optimal input  sequence for
system (\ref{eq:sys}) with cost (\ref{eq:J}).  Hence $\V[d(x)](x)=\J[d(x)](x,\bm{u}^*_{\gamma,d(x)}(x))$.

We define $\Y[d(x)](x):=\V[d(x)](x)+W(x)$ for any $x\in\mathbb{R}^n$, where $W$  comes from Assumption \ref{Aone}. Since item (i) of Theorem 1 in \cite{granzotto2019} is valid for all $d\in\mathbb{Z}_{>0}$,  item (i) of Proposition \ref{YLyapunovProp} holds.
Now, consider the sequence  $\hat{\bm{u}}:=[\us[1],\us[2],\ldots,\us[\bar d-1],\bar{\bm{u}}]$ where
$\bar{\bm{u}}:=\bm{u}^*_{\gamma,\infty}(       \phi({\bar d}       ,      x       ,       \bm{u}^*_{\gamma,d(x)}(x)|_{\bar d})       )$,
$\bm{u}^*_{\gamma,d(x)}(x)|_{\bar d}=[\us[0],\ldots,\us[\bar d-1]]$ and  $\phi$ denotes the  solution of system (\ref{eq:sys}). The sequence
$\hat{\bm{u}}$ consists of the first $\bar d$ elements of  $\bm{u}^*_{\gamma,d(x)}(x)$ after $\us[0]$, followed by an optimal input sequence
of infinite length  at  state $\phi(\bar d,x,\ustar(x)|_{\bar d})$. Note that such sequence only exists and is well defined
if $d(x)\geq \bar d$, which is the case here, and that the sequence $\bar{\bm{u}}$
exists and minimizes  $\J[\infty](\phi(\bar d,x,\ustar(x)|_{\bar d}),\bar{\bm{u}}_{\bar d})$, which is guaranteed to exist by virtue of SA.  From the definition  of cost $\J[d(\cdot)]$
in (\ref{eq:J}) and $\V[d(v)](v)$ in (\ref{eq:V}),
$\V[d(v)](v) \leq
   \J(v,\hat{\bm{u}}) = \J[\bar d-1](v,\hat{\bm{u}}|_{\bar d})+\gamma^{\bar d}\J[\infty](\phi(\bar d,v,\hat{\bm{u}}|_{\bar d}),\bar{\bm{u}}).$
Then, by the same manipulations as in  the proof of Theorem 1 in \cite{granzotto2019}, we obtain
\begin{equation}
\V[d(v)](v)  \leq \frac{\V[d(x)](x)-\ell_{\us[0]}(x)}{\gamma}+\gamma^{\bar d}\overline{\alpha}_V(\sigma(\phi(\bar d,v,\hat{\bm{u}}|_{\bar d}))). \label{eq:propeqfinal}
\end{equation}
We now bound $\sigma(\phi(\bar d,v,\hat{\bm{u}}|_{\bar d}))$ by following the steps of \cite{granzotto2019}. In particular, we have that $\Y[d(x)-(k+1)](\phi(k+1,x,\bm{u}^*_{\gamma,d(x)}(x)|_{k+1}))\leq \left(\frac{\mathbb{I}-\alpha_W\circ\overline{\alpha}_Y^{-1}
}{\gamma}\right)(\Y[d(x)-k](\phi(k,x,\bm{u}^*_{\gamma,d(x)}(x)|_{k}))$ for all $k \in \{1, \ldots,\allowbreak d(x)-1\}$, which yields by iteration, $\Y[d(x)-(k+1)](\phi(k+1,x,\bm{u}^*_{\gamma,d(x)}(x)|_{k+1}))\leq \left(\frac{\mathbb{I}-\alpha_W\circ\overline{\alpha}_Y^{-1}
}{\gamma}\right)^{(k)}(\Y[d(x)](\phi(0,x,\varnothing)))$.  Since  $\sigma(\phi(\bar d,v,\hat{\bm{u}}|_{\bar d}))=\allowbreak
\sigma(\phi(\bar d+1,x,\bm{u}^*_{\gamma,d(x)}(x)|_{\bar d+1}))\leq\allowbreak \alpha_W^{-1}(\Y[d(x)-(\bar d+1)](\phi(\bar d+1,x,\bm{u}^*_{\gamma,d(x)}(x)|_{\bar d+1})))$ per item (i) of Proposition \ref{YLyapunovProp}, it follows by fixing $k=\bar d$, that $\sigma(\phi(\bar d,v,\hat{\bm{u}}|_{\bar d}))\leq\alpha_W^{-1}\left(\left(\frac{\mathbb{I}-\alpha_W\circ\overline{\alpha}_Y^{-1}
}{\gamma}\right)^{(\bar d)}(\Y[d(x)](x))\right)$. The desired result follows by applying the obtained upper-bound to (\ref{eq:propeqfinal}) and noting that $\frac{\V[d(x)](x)}{\gamma}=\V[d(x)](x)+\frac{(1-\gamma)\V[d(x)](x)}{\gamma}$. The case where $\gamma=1$ is similarly treated as \cite{granzotto2019}.
\end{IEEEproof}

We can now finalize the proof of Theorem \ref{algostab}. Since budget
\(B\geq\frac{M^{\bar d+2}-1}{M-1}\), it follows that \(d(x)> \bar d\) by
Proposition \ref{dfiniteoptimum}. By following the proof of Theorem 2 in
\cite{granzotto2019}, the desired result is derived.

\bibliographystyle{plain}

\bibliography{IEEEabrv,OPstability}

\begin{thebibliography}{10}

\bibitem{7588142}
D.~{Antunes} and W.~P. M.~H. {Heemels}.
\newblock Linear quadratic regulation of switched systems using informed
  policies.
\newblock {\em IEEE Transactions on Automatic Control}, 62(6):2675--2688, 2017.

\bibitem{BUSONIU2017297}
L.~Buşoniu, J.~Daafouz, M.~C. Bragagnolo, and I.-C. Morărescu.
\newblock Planning for optimal control and performance certification in
  nonlinear systems with controlled or uncontrolled switches.
\newblock {\em Automatica}, 78:297 -- 308, 2017.

\bibitem{busoniunetwork}
L.~{Buşoniu}, R.~{Postoyan}, and J.~{Daafouz}.
\newblock Near-optimal strategies for nonlinear and uncertain networked control
  systems.
\newblock {\em IEEE Transactions on Automatic Control}, 61(8):2124--2139, 2016.

\bibitem{DEAECTO20181}
G.~S. Deaecto and J.~C. Geromel.
\newblock Stability and performance of discrete-time switched linear systems.
\newblock {\em Systems \& Control Letters}, 118:1 -- 7, 2018.

\bibitem{dellarossatree}
F.~{Della Rossa} and F.~{Dercole}.
\newblock Tree-based algorithms for the stability of discrete-time switched
  linear systems under arbitrary and constrained switching.
\newblock {\em IEEE Transactions on Automatic Control}, 2018.

\bibitem{granzotto2019}
M.~Granzotto, R.~Postoyan, L.~Buşoniu, D.~Nešić, and J.~Daafouz.
\newblock Finite-horizon discounted optimal control: stability and performance.
\newblock In {\em submitted for journal publication}.
\newblock Available for reviewers at \url{https://tinyurl.com/y4y2njr5}.

\bibitem{granzotto2018}
M.~Granzotto, R.~Postoyan, L.~Buşoniu, D.~Nešić, and J.~Daafouz.
\newblock Stability analysis of discrete-time finite-horizon optimal control
  with discounted cost.
\newblock In {\em IEEE Conference on Decision and Control}, Miami, USA, 2018.

\bibitem{grimm2005}
G.~Grimm, M.~J. Messina, S.~E. Tuna, and A.~R. Teel.
\newblock Model predictive control: for want of a local control {Lyapunov}
  function, all is not lost.
\newblock {\em IEEE Transactions on Automatic Control}, 50(5):546--558, 2005.

\bibitem{gruneperformance}
L.~{Gr\"une} and A.~Rantzer.
\newblock On the infinite horizon performance of receding horizon controllers.
\newblock {\em IEEE Transactions on Automatic Control}, 53(9):2100--2111, 2008.

\bibitem{hren2008}
J.-F. Hren and R.~Munos.
\newblock Optimistic planning of deterministic systems.
\newblock In {\em European Workshop on Reinforcement Learning}, pages 151--164,
  Villeneuve d'Ascq, France, 2008. Springer.

\bibitem{keerthi1985}
S.~Keerthi and E.~Gilbert.
\newblock An existence theorem for discrete-time infinite-horizon optimal
  control problems.
\newblock {\em IEEE Transactions on Automatic Control}, 30(9):907--909, 1985.

\bibitem{1039806}
B.~{Lincoln} and B.~{Bernhardsson}.
\newblock {LQR} optimization of linear system switching.
\newblock {\em IEEE Transactions on Automatic Control}, 47(10):1701--1705,
  2002.

\bibitem{munos2014bandits}
R.~Munos et~al.
\newblock From bandits to {Monte-Carlo} tree search: The optimistic principle
  applied to optimization and planning.
\newblock {\em Foundations and Trends{\textregistered} in Machine Learning},
  7(1):1--129, 2014.

\bibitem{romain2016}
R.~Postoyan, L.~Buşoniu, D.~Nešić, and J.~Daafouz.
\newblock Stability analysis of discrete-time infinite-horizon optimal control
  with discounted cost.
\newblock {\em IEEE Transactions on Automatic Control}, 62(6):2736--2749, 2017.

\bibitem{1687497}
A.~{Rantzer}.
\newblock Relaxed dynamic programming in switching systems.
\newblock {\em IEEE Proceedings - Control Theory and Applications},
  153(5):567--574, 2006.

\bibitem{7963352}
J.~B. {Rejeb}, L.~{Buşoniu}, I.~{Morărescu}, and J.~{Daafouz}.
\newblock Near-optimal control of nonlinear switched systems with
  non-cooperative switching rules.
\newblock In {\em American Control Conference (ACC)}, pages 2648--2653,
  Seattle, WA, USA, 2017.

\bibitem{6656817}
P.~{Riedinger}.
\newblock A switched {LQ} regulator design in continuous time.
\newblock {\em IEEE Transactions on Automatic Control}, 59(5):1322--1328, 2014.

\bibitem{4982639}
M.~{Rinehart}, M.~{Dahleh}, and I.~{Kolmanovsky}.
\newblock Value iteration for (switched) homogeneous systems.
\newblock {\em IEEE Transactions on Automatic Control}, 54(6):1290--1294, 2009.

\bibitem{Antsaklis2003}
X.~Xu and P.~J. Antsaklis.
\newblock Results and perspectives on computational methods for optimal control
  of switched systems.
\newblock In O.~Maler and A.~Pnueli, editors, {\em Hybrid Systems: Computation
  and Control}, pages 540--555, Berlin, Heidelberg, 2003. Springer.

\bibitem{5288566}
W.~{Zhang}, J.~{Hu}, and A.~{Abate}.
\newblock On the value functions of the discrete-time switched {LQR} problem.
\newblock {\em IEEE Transactions on Automatic Control}, 54(11):2669--2674,
  2009.

\bibitem{abateDSLQR}
W.~{Zhang}, J.~{Hu}, and A.~{Abate}.
\newblock Infinite-horizon switched {LQR} problems in discrete time: A
  suboptimal algorithm with performance analysis.
\newblock {\em IEEE Transactions on Automatic Control}, 57(7):1815--1821, 2012.

\bibitem{Zhu2015}
F.~Zhu and P.~J. Antsaklis.
\newblock Optimal control of hybrid switched systems: A brief survey.
\newblock {\em Discrete Event Dynamic Systems}, 25(3):345--364, Sep 2015.

\end{thebibliography}

\end{document}